\documentclass{article}
\usepackage{amssymb}
\usepackage{amsfonts}
\usepackage{hyperref}
\usepackage{amsmath}
\usepackage{geometry}

\newtheorem{theorem}{Theorem}

\newtheorem{conjecture}[theorem]{Conjecture}
\newtheorem{corollary}[theorem]{Corollary}

\newtheorem{definition}[theorem]{Definition}

\newtheorem{lemma}[theorem]{Lemma}

\newtheorem{proposition}[theorem]{Proposition}
\newtheorem{remark}[theorem]{Remark}

\numberwithin{theorem}{section}
\input{tcilatex}
\begin{document}

\title{An Explicit Presentation of the Grothendieck Ring of Finitely
Generated $\mathbb{F}_{q}[SL_{2}(\mathbb{F}_{q})]$-Modules}
\author{Davide A. Reduzzi \\
University of California at Los Angeles\\
{\scriptsize devredu83@math.ucla.edu}}
\maketitle

\begin{abstract}
Let $p$ be a prime and $q=p^{g}$. We show that the Grothendieck ring of
finitely generated $\mathbb{F}_{q}[SL_{2}(\mathbb{F}_{q})]$-modules is
naturally isomorphic to the quotient of the polynomial algebra $%
\mathbb{Z}
\lbrack x]$ by the ideal generated by $f^{[g]}(x)-x$, where $%
f(x)=\tsum\nolimits_{j=0}^{\left\lfloor p/2\right\rfloor }\left( -1\right)
^{j}\tfrac{p}{p-j}\tbinom{p-j}{j}x^{p-2j}$, and the superscript $[g]$
denotes $g$-fold composition of polynomials. We conjecture that a similar
result holds for simply connected semisimple algebraic groups defined and
split over a finite field.
\end{abstract}

\section{Introduction}

In \cite{Se01}, J-P. Serre discovered a puzzling identity involving
characteristic $p$ symmetric powers representations of the group $GL_{2}(%
\mathbb{F}_{q})$, viewed as elements of the Grothendieck ring $K_{0}\left(
GL_{2}(\mathbb{F}_{q})\right) $ of finitely generated $\mathbb{F}_{q}[GL_{2}(%
\mathbb{F}_{q})]$-modules.

More precisely, fix a rational prime $p$, a positive integer $g$, and set $%
q=p^{g}$. Denote by $\mathbb{F}_{q}$ a field with $q$ elements and by $%
\mathfrak{G}$ the group $SL_{2}(\mathbb{F}_{q})\subset GL_{2}(\mathbb{F}%
_{q}) $. For any non-negative integer $k$, denote by $\mathfrak{M}_{k}$ the $%
(k+1)$-dimensional representation $\limfunc{Sym}\nolimits^{k}\mathbb{F}%
_{q}^{2}$ of $\mathfrak{G}$. Motivated by the computation of the
Euler-Poincar\'{e} characteristic of the twisted sheaf $\mathcal{O}\left(
k\right) $ on $\mathbb{P}_{\mathbb{F}_{q}}^{1}$, in \cite{Se01} Serre
extended the definition of the modules $\mathfrak{M}_{k}$'s for negative
values of $k$, and showed that for any integer $k$ the following relation
holds in the ring $K_{0}\left( \mathfrak{G}\right) $:%
\begin{equation}
\mathfrak{M}_{k}-\mathfrak{M}_{k-(q+1)}=\mathfrak{M}_{k-\left( q-1\right) }-%
\mathfrak{M}_{k-2q}.  \tag{$\Sigma $}
\end{equation}

The dimensional shiftings by $q+1$ and $q-1$ occurring in Serre's relation
can be obtained by applying opportune intertwining operators $\Theta _{q}$
and $D$ to the symmetric powers modules. This has been exploited in \cite{Re}%
\ for the study of cohomological weight shiftings for elliptic modular forms
modulo $p$.

Motivated by generalizations of the above considerations to Hilbert modular
forms, families of generalized $\Theta _{q}$ and $D$ operators are defined
in \cite{Re2}, and the following identity in $K_{0}\left( \mathfrak{G}%
\right) $ is proved for any integers $k,h$ and $i$:%
\begin{equation}
\mathfrak{M}_{k}^{[i]}\mathfrak{M}_{h}^{[i+1]}-\mathfrak{M}_{k-p}^{[i]}%
\mathfrak{M}_{h-1}^{[i+1]}=\mathfrak{M}_{k-p}^{[i]}\mathfrak{M}%
_{h+1}^{[i+1]}-\mathfrak{M}_{k-2p}^{[i]}\mathfrak{M}_{h}^{[i+1]}. 
\tag{$\Phi $}
\end{equation}

\noindent Here the superscript $[i]$ denote the $i$th Frobenius twisting on
the corresponding virtual representation.

Using Glover's product identity, one sees that $\left( \Phi \right) $ is
equivalent to $\left( \Sigma \right) $ in case $g=1$, but it is stronger for 
$g>1.$

In this paper we apply formula $(\Phi )$ to determine an explicit
presentation of the Grothendieck ring $K_{0}\left( \mathfrak{G}\right) $. We
treat the case of $\mathfrak{G}=SL_{2}(\mathbb{F}_{q})$ instead of $GL_{2}(%
\mathbb{F}_{q})$, so we will not need to consider determinant twists that
would make the set of relations more complicated; following the same methods
we describe below, one could easily work with $GL_{2}(\mathbb{F}_{q})$
instead.

Our main result is the following (cf. Theorem \ref{presentation}):

\begin{theorem}
Denote by $\mathfrak{X}$ the standard representation of $\mathfrak{G}$ on $%
\mathbb{F}_{q}^{2}$ and let $\mathfrak{x}$ be an indeterminate over $%
\mathbb{Z}
$. The assignment $\mathfrak{X\longmapsto x}$\ induces an isomorphism of
rings:%
\begin{equation*}
K_{0}\left( \mathfrak{G}\right) \simeq \frac{%
\mathbb{Z}
\lbrack \mathfrak{x}]}{\left( \mathfrak{f}^{[g]}(\mathfrak{x})-\mathfrak{x}%
\right) 
\mathbb{Z}
\lbrack \mathfrak{x}]},
\end{equation*}

\noindent where $\mathfrak{f}^{[g]}(\mathfrak{x})=\left( \mathfrak{f}\circ
...\circ \mathfrak{f}\right) \left( \mathfrak{x}\right) $ is the polynomial
of $%
\mathbb{Z}
\lbrack \mathfrak{x}]$ having degree $p^{g}$ obtained by $g$-fold
composition of the monic degree $p$ polynomial:%
\begin{equation*}
\mathfrak{f}(\mathfrak{x})=\dsum\nolimits_{j=0}^{\left\lfloor
p/2\right\rfloor }\left( -1\right) ^{j}\dfrac{p}{p-j}\dbinom{p-j}{j}%
\mathfrak{x}^{p-2j}.
\end{equation*}
\end{theorem}

Proposition \ref{gth compo} gives an explicit closed formula for $\mathfrak{f%
}^{[g]}(\mathfrak{x})$. Notice that, since $\mathfrak{f}(\mathfrak{x})\equiv 
\mathfrak{x}^{p}(\func{mod}p%
\mathbb{Z}
\lbrack \mathfrak{x}])$, the structure of the generic and special fibers of
the ring $K_{0}\left( \mathfrak{G}\right) \otimes _{%
\mathbb{Z}
}%
\mathbb{Z}
_{p}$ are easily determined (Corollary \ref{cor_presentation}). On the other
side, the arithmetic properties of the polynomial $\mathfrak{f}(\mathfrak{x}%
) $ over $%
\mathbb{Q}
$ seem to be more complicated.

In the last paragraph of the paper we prove the following fact (Proposition %
\ref{tensorp}): assume $\mathbb{G}$ is a simply connected, semisimple
algebraic group defined and split over $\mathbb{F}_{q}$. If $\mathfrak{M}$
is an $\mathbb{F}_{q}[\mathbb{G]}$-rational module of finite $\mathbb{F}_{q}$%
-dimension, then the multiplicity of an irreducible $\mathbb{F}_{q}[\mathbb{%
G]}$-rational module $\mathfrak{V}$ as a Jordan-H\"{o}lder constituent of $%
\mathfrak{M}^{[i]}$ is congruent modulo $p$ to the multiplicity of $%
\mathfrak{V}$ as a Jordan-H\"{o}lder constituent of $\mathfrak{M}^{\otimes
p^{i}}$, for any positive integer $i$.

Motivated by this result, we are led to conjecture that the Grothendieck
ring of a Chevalley group arising from a rank $\ell $ algebraic group $%
\mathbb{G}$ as above is isomorphic to the algebra%
\begin{equation*}
\frac{%
\mathbb{Z}
\lbrack \mathfrak{x}_{1},...,\mathfrak{x}_{\ell }]}{\left( \mathfrak{f}%
_{1}^{[g]}(\mathfrak{x}_{1})-\mathfrak{x}_{1}\mathfrak{,...,f}_{\ell }^{[g]}(%
\mathfrak{x}_{\ell })-\mathfrak{x}_{\ell }\right) 
\mathbb{Z}
\lbrack \mathfrak{x}_{1},...,\mathfrak{x}_{\ell }]},
\end{equation*}

\noindent where for any $i$, $\mathfrak{f}_{i}^{[g]}(\mathfrak{x}_{i})$ is
the $g$-fold composition of the degree $p$ monic polynomial $\mathfrak{f}%
_{i}(\mathfrak{x}_{i})\in 
\mathbb{Z}
\lbrack \mathfrak{x}_{i}]$. We conjecture that $\mathfrak{f}_{i}(\mathfrak{x}%
_{i})\equiv \mathfrak{x}_{i}^{p}(\func{mod}p%
\mathbb{Z}
\lbrack \mathfrak{x}_{i}])$ for any value of $i$. Some of the evidence for
this conjecture is presented at the end of paragraph 4.

\bigskip

\textbf{Conventions} All the group representations in this paper are left
representations on a module of finite length over a fixed ring. If $R$ is an
algebra over a ring $A$, and $\mathcal{S}$ is a subset of $R$, the symbol $A[%
\mathcal{S}]$ denotes the $A$-subalgebra of $R$ generated by $\mathcal{S}$.

\bigskip

\section{Weight shiftings identities in $K_{0}\left( \mathfrak{G}\right) $}

Fix a rational prime $p$, a positive integer $g$, and set $q=p^{g}$. Denote
by $\mathbb{F}_{q}$ a finite field with $q$ elements and fix an algebraic
closure $\overline{\mathbb{F}}_{q}$ of $\mathbb{F}_{q}$; let $\sigma \in 
\limfunc{Gal}\left( \mathbb{F}_{q}/\mathbb{F}_{p}\right) $ be the arithmetic
absolute Frobenius element. Denote by $\mathfrak{G}$ the group $SL_{2}(%
\mathbb{F}_{q})$.

For any $i\in 
\mathbb{Z}
$, the Frobenius power $\sigma ^{i}$ induces a function $\mathfrak{G}%
\longrightarrow \mathfrak{G}$ obtained by applying $\sigma ^{i}$ to each
entry of the matrices in $\mathfrak{G}$: composing this map with the action
of $\mathfrak{G}$ on a given $\mathbb{F}_{q}[\mathfrak{G}]$-module $%
\mathfrak{M}$ gives to the latter a new structure of $\mathfrak{G}$-module,
that is denoted $\mathfrak{M}^{[i]}$ and called the $i$th Frobenius twist of 
$\mathfrak{M}.$

If $f:\mathfrak{M}\longrightarrow \mathfrak{N}$ is a homomorphism of $%
\mathbb{F}_{q}[\mathfrak{G}]$-modules and $i\in 
\mathbb{Z}
$, we denote by $f^{[i]}:\mathfrak{M}^{[i]}\longrightarrow \mathfrak{N}%
^{[i]} $ the $\mathfrak{G}$-homomorphism defined by $f^{[i]}(x)=f(x)$ for
all $x\in \mathfrak{M}^{[i]}$.

Let $\mathfrak{X}$ denote the standard representation of $\mathfrak{G}$ on $%
\mathbb{F}_{q}^{2}$ and, for any positive integer $k$, define 
\begin{equation*}
\mathfrak{M}_{k}=\limfunc{Sym}\nolimits^{k}\mathfrak{X}
\end{equation*}%
to be the $k$th symmetric power of $\mathfrak{X}$, so that in particular $%
\mathfrak{X=M}_{1}$. Let $\mathfrak{M}_{0}$ be the trivial representation of 
$\mathfrak{G}$ on $\mathbb{F}_{q}$.

Observe that we can identify $\mathfrak{M}_{k}$ with the $\mathbb{F}_{q}$%
-vector space of homogeneous degree $k$ polynomials over $\mathbb{F}_{q}$ in
two variables $X$ and $Y$, endowed with the action of $\mathfrak{G}$ induced
by:

\begin{center}
\begin{equation*}
\left( 
\begin{array}{cc}
a & b \\ 
c & d%
\end{array}%
\right) \cdot X=aX+cY,\ \ \ \left( 
\begin{array}{cc}
a & b \\ 
c & d%
\end{array}%
\right) \cdot Y=bX+dY.
\end{equation*}
\end{center}

As a consequence of Steinberg's restriction and tensor product theorems (%
\cite{St}) we have:

\begin{proposition}
\label{stei}All and only the irreducible representations of $\mathfrak{G}$
over $\mathbb{F}_{q}$ are of the form:%
\begin{equation*}
\dbigotimes\nolimits_{i=0}^{g-1}\mathfrak{M}_{k_{i}}^{[i]},
\end{equation*}
where $k_{0},...,k_{g-1}$ are integers such that $0\leq k_{i}\leq p-1$ for $%
i=0,...,g-1$, and the tensor products are over $\mathbb{F}_{q}$.
Furthermore, the above representations are pairwise non-isomorphic.
\end{proposition}

Denote by $K_{0}(\mathfrak{G})$ the Grothendieck group of finitely generated 
$\mathbb{F}_{q}[\mathfrak{G}]$-modules. $K_{0}(\mathfrak{G})$ is isomorphic
to the free abelian group generated by the isomorphism classes of
irreducible representations of $\mathfrak{G}$ over $\mathbb{F}_{q}$, so that
it has rank $q$ over $%
\mathbb{Z}
$. If $\mathfrak{M}$ is an $\mathbb{F}_{q}[\mathfrak{G}]$-module, we will
denote by the same symbol $\mathfrak{M}$ its class in $K_{0}(\mathfrak{G})$,
if no confusion arises.

Tensor product over $\mathbb{F}_{q}$ induces on $K_{0}(\mathfrak{G})$ a
structure of commutative ring with identity; we denote the product in $K_{0}(%
\mathfrak{G})$ by $\cdot $ or by simple juxtaposition. All the tensor
products we will consider in the sequel are over $\mathbb{F}_{q}$, unless
otherwise specified.

\bigskip

We can extend the definition of the virtual representation $\mathfrak{M}_{k}$
for $k<0$ in a way that is coherent with Brauer character computations. In 
\cite{Se01}, the determination the Euler-Poincar\'{e} characteristic of the
twisted sheaf $\mathcal{O}\left( k\right) $ on $\mathbb{P}_{\mathbb{F}%
_{q}}^{1}$ suggests the following:

\begin{definition}
Let $k$ be a negative integer. Define the element $\mathfrak{M}_{k}$ of the
Grothendieck group $K_{0}\left( \mathfrak{G}\right) $ of $\mathfrak{G}$ over 
$\mathbb{F}_{q}$ by:%
\begin{equation*}
\mathfrak{M}_{k}=\left\{ 
\begin{array}{cc}
0 & \text{if }k=-1 \\ 
\mathfrak{M}_{-k-2} & \text{if }k\leq -2.%
\end{array}%
\right.
\end{equation*}
\end{definition}

The following result summarizes some non-trivial identities that hold in the
ring $K_{0}\left( \mathfrak{G}\right) $:

\begin{theorem}
Let $k$ and $h$ be any integers. The following identities hold in $%
K_{0}\left( \mathfrak{G}\right) $:%
\begin{equation}
\mathfrak{M}_{k}=-\mathfrak{M}_{-k-2}  \tag{$\Delta _{g}$}
\end{equation}%
\begin{equation}
\mathfrak{M}_{k}-\mathfrak{M}_{k-(q+1)}=\mathfrak{M}_{k-\left( q-1\right) }-%
\mathfrak{M}_{k-2q}  \tag{$\Sigma _{g}$}
\end{equation}%
\begin{equation}
\mathfrak{M}_{k}\mathfrak{M}_{h}=\mathfrak{M}_{k+h}+\mathfrak{M}_{k-1}%
\mathfrak{M}_{h-1}  \tag{$\Pi _{g}$}
\end{equation}%
\begin{equation}
\mathfrak{M}_{k}=\mathfrak{M}_{k-p}\mathfrak{M}_{1}^{[1]}-\mathfrak{M}%
_{k-2p}.  \tag{$\Phi _{g}$}
\end{equation}
\end{theorem}

\textbf{Proof }Formulae $(\Delta _{g})$ and $(\Sigma _{g})$ are proved in 
\cite{Se01} via a Brauer characters computation. Formula $(\Pi _{g})$ comes
from an exact sequence of $\mathfrak{G}$-modules constructed in \cite{Glo}.
Formula $(\Phi _{g})$ is proved in section 3 of \cite{Re2}. $\blacksquare $

\bigskip

We remark that formula $(\Phi _{g})$ appeared in \cite{Re2} also in the form:%
\begin{equation*}
\mathfrak{M}_{k}^{[i]}\mathfrak{M}_{h}^{[i+1]}-\mathfrak{M}_{k-p}^{[i]}%
\mathfrak{M}_{h-1}^{[i+1]}=\mathfrak{M}_{k-p}^{[i]}\mathfrak{M}%
_{h+1}^{[i+1]}-\mathfrak{M}_{k-2p}^{[i]}\mathfrak{M}_{h}^{[i+1]},
\end{equation*}

\noindent where $k,h$ and $i$ are any integers.

The product formula $(\Pi _{1})$ implies that $(\Phi _{1})$ and $\left(
\Sigma _{1}\right) $ are equivalent. If $g>1$, $(\Phi _{g})$ cannot be
deduced from $(\Sigma _{g})$ and $(\Pi _{g})$: the proof of this fact,
contained in \cite{Re2}, is indirect and throughout the paper the knowledge
of Serre's relation $(\Sigma _{g})$ will allow sometimes to bypass long
computations involving Frobenius twists, when $g>1$.

In \cite{Re2} it is also proved that for $g\geq 1$, we can use the relations 
$(\Delta _{g}),(\Phi _{g}),(\Pi _{g})$ to explicitly compute the Jordan-H%
\"{o}lder factors of any virtual representations of the form $%
\tprod\nolimits_{i=0}^{g-1}\mathfrak{M}_{k_{i}}^{[i]}$, where $%
k_{0},...,k_{g-1}\in 
\mathbb{Z}
$.

\section{Presentation of $K_{0}\left( \mathfrak{G}\right) $}

We keep the notation introduced in the previous paragraph.

\begin{lemma}
\label{1generator}The ring $K_{0}\left( \mathfrak{G}\right) $ is generated
by $\mathfrak{X}$ as a $%
\mathbb{Z}
$-algebra.
\end{lemma}

\textbf{Proof }By Proposition \ref{stei}, $K_{0}\left( \mathfrak{G}\right) $
is freely\ generated as a $%
\mathbb{Z}
$-module by the $q$ elements $\tprod\nolimits_{i=0}^{g-1}\mathfrak{M}%
_{k_{i}}^{[i]}$, where $0\leq k_{i}\leq p-1$ for any $i$. It is therefore
enough to show that for all integers $i,k$ such that $0\leq i\leq g-1$ and $%
0\leq k\leq p-1$ we have $\mathfrak{M}_{k}^{[i]}\in 
\mathbb{Z}
\lbrack \mathfrak{X}]$.

Applying $(\Pi _{g})$ we obtain the recursive relations:%
\begin{equation}
\mathfrak{M}_{2}=\mathfrak{X}^{2}-1,\text{ }\mathfrak{M}_{n}=\mathfrak{%
X\cdot M}_{n-1}-\mathfrak{M}_{n-2}\text{ \ (}n>2\text{),}  \label{m}
\end{equation}
so that $\mathfrak{M}_{k}\in 
\mathbb{Z}
\lbrack \mathfrak{X}]$ for all $k\geq 0$. Twisting (\ref{m})\ by powers of
Frobenius, we obtain:

\begin{equation*}
\mathfrak{M}_{2}^{[i]}=\left( \mathfrak{X}^{[i]}\right) ^{2}-1,\text{ }%
\mathfrak{M}_{n}^{[i]}=\mathfrak{X}^{[i]}\cdot \mathfrak{M}_{n-1}^{[i]}-%
\mathfrak{M}_{n-2}^{[i]}\text{ \ (}n>2\text{),}
\end{equation*}

\noindent for all $0\leq i\leq g-1$, so that $\mathfrak{M}_{k}^{[i]}\in 
\mathbb{Z}
\lbrack \mathfrak{X,X}^{[1]},...,\mathfrak{X}^{[g-1]}]$ for all $k\geq 0$
and: 
\begin{equation*}
K_{0}\left( \mathfrak{G}\right) =%
\mathbb{Z}
\lbrack \mathfrak{X,X}^{[1]},...,\mathfrak{X}^{[g-1]}].
\end{equation*}
By $(\Phi _{g})$, we have $\mathfrak{M}_{p}=\mathfrak{M}_{1}^{[1]}-\mathfrak{%
M}_{-p}$, and applying $(\Delta _{g})$ we obtain $\mathfrak{X}^{[1]}=%
\mathfrak{M}_{p}-\mathfrak{M}_{p-2}$, so that $\mathfrak{X}^{[1]}\in 
\mathbb{Z}
\lbrack \mathfrak{X}]$, as $\mathfrak{M}_{k}\in 
\mathbb{Z}
\lbrack \mathfrak{X}]$ for all $k\geq 0$. We also obtain that, for any $%
0\leq i\leq g-1$, we have:%
\begin{equation}
\mathfrak{X}^{[i+1]}=\mathfrak{M}_{p}^{[i]}-\mathfrak{M}_{p-2}^{[i]},
\label{ff}
\end{equation}%
and we conclude $\mathfrak{X}^{[1]},...,\mathfrak{X}^{[g-1]}\in 
\mathbb{Z}
\lbrack \mathfrak{X}]$, implying $K_{0}\left( \mathfrak{G}\right) =%
\mathbb{Z}
\lbrack \mathfrak{X}]$. $\blacksquare $

\bigskip

Let $\mathfrak{x}$ be an indeterminate over $%
\mathbb{Z}
$ and define the following two families of polynomials of $%
\mathbb{Z}
\lbrack \mathfrak{x}]$:%
\begin{eqnarray*}
&&\left\{ 
\begin{array}{l}
\mathfrak{m}_{0}(\mathfrak{x})=1 \\ 
\mathfrak{m}_{1}(\mathfrak{x})=\mathfrak{x} \\ 
\mathfrak{m}_{2}(\mathfrak{x})=\mathfrak{x}^{2}-1 \\ 
\mathfrak{m}_{n}(\mathfrak{x})=\mathfrak{x\cdot m}_{n-1}(\mathfrak{x})-%
\mathfrak{m}_{n-2}(\mathfrak{x})\ \text{\ \ }(n>2);%
\end{array}%
\right. \\
&& \\
&&\left\{ 
\begin{array}{l}
\mathfrak{f}^{[0]}(\mathfrak{x})=\mathfrak{x} \\ 
\mathfrak{f}^{[1]}(\mathfrak{x})=\mathfrak{m}_{p}(\mathfrak{x})-\mathfrak{m}%
_{p-2}(\mathfrak{x}) \\ 
\mathfrak{f}^{[i]}(\mathfrak{x})=\text{ }\left( \mathfrak{f}^{[1]}\circ 
\mathfrak{f}^{[i-1]}\right) (\mathfrak{x})\ =\mathfrak{m}_{p}(\mathfrak{f}%
^{[i-1]}(\mathfrak{x}))-\mathfrak{m}_{p-2}(\mathfrak{f}^{[i-1]}(\mathfrak{x}%
))\text{\ \ \ }(i>1).%
\end{array}%
\right.
\end{eqnarray*}

Observe that for any non-negative integer $n$, $\mathfrak{m}_{n}(\mathfrak{x}%
)$ is a monic polynomial of degree $n$, so that for any non-negative integer 
$i$, $\mathfrak{f}^{[i]}(\mathfrak{x})$ is a monic polynomial of degree $%
p^{i}$.

\begin{lemma}
For any non-negative integer $i$, we have $\mathfrak{f}^{[i]}(\mathfrak{X})=%
\mathfrak{X}^{[i]}$ in $K_{0}\left( \mathfrak{G}\right) $.
\end{lemma}

\textbf{Proof }Notice first that, by definition of $\mathfrak{m}_{n}(%
\mathfrak{x})$ and by formula (\ref{m}), one has:%
\begin{equation}
\mathfrak{m}_{n}(\mathfrak{X})=\mathfrak{M}_{n}  \label{meaningful1}
\end{equation}

\noindent in $K_{0}\left( \mathfrak{G}\right) $ ($n\geq 0$). To prove the
lemma, we use induction on $i$. If $i=0$, the statement is clear; if $i=1$
it follows from formulae (\ref{meaningful1}) and (\ref{ff}). Assume $i\geq 1$
fixed and suppose $\mathfrak{f}^{[i]}(\mathfrak{X})=\mathfrak{X}^{[i]}$. We
have:%
\begin{eqnarray*}
\mathfrak{f}^{[i+1]}(\mathfrak{X}) &=&\mathfrak{m}_{p}(\mathfrak{f}^{[i]}(%
\mathfrak{X}))-\mathfrak{m}_{p-2}(\mathfrak{f}^{[i]}(\mathfrak{X})) \\
&=&\mathfrak{m}_{p}(\mathfrak{X}^{[i]})-\mathfrak{m}_{p-2}(\mathfrak{X}%
^{[i]}).
\end{eqnarray*}

\noindent Observe that Frobenius twists do not act on the coefficients of
virtual representations in $K_{0}\left( \mathfrak{G}\right) $, so that the
last term above is equal to $\mathfrak{m}_{p}(\mathfrak{X})^{[i]}-\mathfrak{m%
}_{p-2}(\mathfrak{X})^{[i]}$. By formula (\ref{meaningful1}), the latter is $%
\mathfrak{M}_{p}^{[i]}-\mathfrak{M}_{p-2}^{[i]}$. By formula (\ref{ff}),
this is $\mathfrak{M}_{p}^{[i]}-\mathfrak{M}_{p-2}^{[i]}=\mathfrak{X}%
^{[i+1]} $. $\blacksquare $

\bigskip

\begin{proposition}
There is an isomorphism of rings:%
\begin{equation*}
\frac{%
\mathbb{Z}
\lbrack \mathfrak{x}]}{\left( \mathfrak{f}^{[g]}(\mathfrak{x})-\mathfrak{x}%
\right) 
\mathbb{Z}
\lbrack \mathfrak{x}]}\simeq K_{0}\left( \mathfrak{G}\right) ,
\end{equation*}

\noindent induced by mapping the indeterminate $\mathfrak{x}$ of the
polynomial ring $%
\mathbb{Z}
\lbrack \mathfrak{x}]$ into the class of the representation $\mathfrak{X}$
of $\mathfrak{G}$.
\end{proposition}

\textbf{Proof }By Proposition \ref{1generator}, the ring homomorphism $%
\mathbb{Z}
\lbrack \mathfrak{x}]\longrightarrow K_{0}\left( \mathfrak{G}\right) $
induced by $\mathfrak{x}\mapsto \mathfrak{X}$ is surjective. Since $%
\mathfrak{X}^{[g]}=\mathfrak{X}$ in $K_{0}\left( \mathfrak{G}\right) $, and
since by the above lemma we have $\mathfrak{f}^{[g]}(\mathfrak{X})=\mathfrak{%
X}^{[g]}$, the above assignment induces an epimorphism 
\begin{equation*}
\pi :\frac{%
\mathbb{Z}
\lbrack \mathfrak{x}]}{\left( \mathfrak{f}^{[g]}(\mathfrak{x})-\mathfrak{x}%
\right) 
\mathbb{Z}
\lbrack \mathfrak{x}]}\longrightarrow K_{0}\left( \mathfrak{G}\right) .
\end{equation*}%
Since $\mathfrak{f}^{[g]}(\mathfrak{x})-\mathfrak{x}$ is a polynomial of
degree $p^{g}$ and since $K_{0}\left( \mathfrak{G}\right) $ is $%
\mathbb{Z}
$-free of rank $p^{g}$, after tensoring with $%
\mathbb{Q}
$ the map $\pi $ defines an isomorphism of $%
\mathbb{Q}
$-vector spaces. This implies that $\ker \pi $ is a finitely generated
torsion $%
\mathbb{Z}
$-submodule of $\frac{%
\mathbb{Z}
\lbrack \mathfrak{x}]}{\left( \mathfrak{f}^{[g]}(\mathfrak{x})-\mathfrak{x}%
\right) 
\mathbb{Z}
\lbrack \mathfrak{x}]}$, and hence it is trivial since $\mathfrak{f}^{[g]}(%
\mathfrak{x})-\mathfrak{x}$ is monic. We conclude that $\pi $ is an
isomorphism of rings. $\blacksquare $

\bigskip

We are now left with determining an explicit formula for the polynomial $%
\mathfrak{f}^{[g]}(\mathfrak{x})\in 
\mathbb{Z}
\lbrack \mathfrak{x}]$.

\begin{lemma}
For any non-negative integer $n$ we have: 
\begin{equation*}
\mathfrak{m}_{n}(\mathfrak{x})=\dsum\nolimits_{j=0}^{\left\lfloor
n/2\right\rfloor }\left( -1\right) ^{j}\dbinom{n-j}{j}\mathfrak{x}^{n-2j}.
\end{equation*}

\noindent \noindent (Where, for any integer $h$, $\left\lfloor
h\right\rfloor $ denotes the largest integer not greater than $h$).
\end{lemma}

\textbf{Proof }We use induction on $n\geq 0$; denote by $\mathfrak{m}%
_{n}^{\prime }(\mathfrak{x})$ the right hand side of the above formula. We
have $\mathfrak{m}_{0}^{\prime }(\mathfrak{x})=1=\mathfrak{m}_{0}(\mathfrak{x%
})$, $\mathfrak{m}_{1}^{\prime }(\mathfrak{x})=\mathfrak{x}=\mathfrak{m}_{1}(%
\mathfrak{x})$ and $\mathfrak{m}_{2}^{\prime }(\mathfrak{x})=\mathfrak{x}%
^{2}-1=\mathfrak{m}_{2}(\mathfrak{x})$. If $n>2$ we have by induction:

\begin{eqnarray*}
\mathfrak{m}_{n}(\mathfrak{x}) &=&\mathfrak{xm}_{n-1}(\mathfrak{x})-%
\mathfrak{m}_{n-2}(\mathfrak{x})=\mathfrak{xm}_{n-1}^{\prime }(\mathfrak{x})-%
\mathfrak{m}_{n-2}^{\prime }(\mathfrak{x}) \\
&=&\tsum\nolimits_{j=0}^{\left\lfloor (n-1)/2\right\rfloor }\left( -1\right)
^{j}\tbinom{n-1-j}{j}\mathfrak{x}^{n-2j}-\tsum\nolimits_{j=0}^{\left\lfloor
(n-2)/2\right\rfloor }\left( -1\right) ^{j}\tbinom{n-2-j}{j}\mathfrak{x}%
^{n-2(j+1)}.
\end{eqnarray*}

\noindent If $n>2$ is even, $\left\lfloor (n-1)/2\right\rfloor =\left\lfloor
(n-2)/2\right\rfloor =(n/2)-1$ and:

\begin{eqnarray*}
\mathfrak{m}_{n}(\mathfrak{x}) &=&\tsum\nolimits_{j=0}^{(n/2)-1}\left(
-1\right) ^{j}\tbinom{n-1-j}{j}\mathfrak{x}^{n-2j}+\tsum%
\nolimits_{j=1}^{n/2}\left( -1\right) ^{j}\tbinom{n-1-j}{j-1}\mathfrak{x}%
^{n-2j} \\
&=&\mathfrak{x}^{n}+\left( \tsum\nolimits_{j=1}^{(n/2)-1}\left( -1\right)
^{j}\left[ \tbinom{n-1-j}{j}+\tbinom{n-1-j}{j-1}\right] \mathfrak{x}%
^{n-2j}\right) +(-1)^{n/2} \\
&=&\mathfrak{x}^{n}+\left( \tsum\nolimits_{j=1}^{(n/2)-1}\left( -1\right)
^{j}\tbinom{n-j}{j}\mathfrak{x}^{n-2j}\right) +(-1)^{n/2} \\
&=&\tsum\nolimits_{j=0}^{\left\lfloor n/2\right\rfloor }\left( -1\right) ^{j}%
\tbinom{n-j}{j}\mathfrak{x}^{n-2j} \\
&=&\mathfrak{m}_{n}^{\prime }(\mathfrak{x}).
\end{eqnarray*}

\noindent If $n>2$ is odd, $\left\lfloor (n-1)/2\right\rfloor =(n-1)/2$, $%
\left\lfloor (n-2)/2\right\rfloor =(n-3)/2$ and:

\begin{eqnarray*}
\mathfrak{m}_{n}(\mathfrak{x}) &=&\tsum\nolimits_{j=0}^{(n-1)/2}\left(
-1\right) ^{j}\tbinom{n-1-j}{j}\mathfrak{x}^{n-2j}+\tsum%
\nolimits_{j=1}^{(n-1)/2}\left( -1\right) ^{j}\tbinom{n-1-j}{j-1}\mathfrak{x}%
^{n-2j} \\
&=&\mathfrak{x}^{n}+\left( \tsum\nolimits_{j=1}^{(n-1)/2}\left( -1\right)
^{j}\left[ \tbinom{n-1-j}{j}+\tbinom{n-1-j}{j-1}\right] \mathfrak{x}%
^{n-2j}\right) \\
&=&\mathfrak{x}^{n}+\left( \tsum\nolimits_{j=1}^{(n-1)/2}\left( -1\right)
^{j}\tbinom{n-j}{j}\mathfrak{x}^{n-2j}\right) \\
&=&\tsum\nolimits_{j=0}^{\left\lfloor n/2\right\rfloor }\left( -1\right) ^{j}%
\tbinom{n-j}{j}\mathfrak{x}^{n-2j} \\
&=&\mathfrak{m}_{n}^{\prime }(\mathfrak{x}).\text{ \ }\blacksquare
\end{eqnarray*}

\bigskip

\begin{corollary}
\label{n-n2}Let $n\geq 2$ be an integer. Then:%
\begin{equation*}
\mathfrak{m}_{n}(\mathfrak{x})-\mathfrak{m}_{n-2}(\mathfrak{x}%
)=\dsum\nolimits_{j=0}^{\left\lfloor n/2\right\rfloor }\left( -1\right) ^{j}%
\dfrac{n}{n-j}\dbinom{n-j}{j}\mathfrak{x}^{n-2j}.
\end{equation*}
\end{corollary}

\textbf{Proof }This is a computation using the previous lemma. We
distinguish two cases: if $n\geq 2$ is even we have:

\begin{eqnarray*}
\mathfrak{m}_{n}(\mathfrak{x})-\mathfrak{m}_{n-2}(\mathfrak{x})
&=&\tsum\nolimits_{j=0}^{n/2}\left( -1\right) ^{j}\tbinom{n-j}{j}\mathfrak{x}%
^{n-2j}-\tsum\nolimits_{j=0}^{(n/2)-1}\left( -1\right) ^{j}\tbinom{n-2-j}{j}%
\mathfrak{x}^{n-2(j+1)} \\
&=&\tsum\nolimits_{j=0}^{n/2}\left( -1\right) ^{j}\tbinom{n-j}{j}\mathfrak{x}%
^{n-2j}+\tsum\nolimits_{j=1}^{n/2}\left( -1\right) ^{j}\tbinom{n-1-j}{j-1}%
\mathfrak{x}^{n-2j} \\
&=&\mathfrak{x}^{n}+\tsum\nolimits_{j=1}^{n/2}\left( -1\right) ^{j}\left[ 
\tbinom{n-j}{j}+\tbinom{n-1-j}{j-1}\right] \mathfrak{x}^{n-2j} \\
&=&\mathfrak{x}^{n}+\tsum\nolimits_{j=1}^{n/2}\left( -1\right) ^{j}\tfrac{n}{%
n-j}\tbinom{n-j}{j}\mathfrak{x}^{n-2j}.
\end{eqnarray*}

\noindent If $n\geq 3$ is odd we have:

\begin{eqnarray*}
\mathfrak{m}_{n}(\mathfrak{x})-\mathfrak{m}_{n-2}(\mathfrak{x})
&=&\tsum\nolimits_{j=0}^{(n-1)/2}\left( -1\right) ^{j}\tbinom{n-j}{j}%
\mathfrak{x}^{n-2j}-\tsum\nolimits_{j=0}^{(n-3)/2}\left( -1\right) ^{j}%
\tbinom{n-2-j}{j}\mathfrak{x}^{n-2(j+1)} \\
&=&\tsum\nolimits_{j=0}^{(n-1)/2}\left( -1\right) ^{j}\tbinom{n-j}{j}%
\mathfrak{x}^{n-2j}+\tsum\nolimits_{j=1}^{(n-1)/2}\left( -1\right) ^{j}%
\tbinom{n-1-j}{j-1}\mathfrak{x}^{n-2j} \\
&=&\mathfrak{x}^{n}+\tsum\nolimits_{j=1}^{(n-1)/2}\left( -1\right) ^{j}\left[
\tbinom{n-j}{j}+\tbinom{n-1-j}{j-1}\right] \mathfrak{x}^{n-2j} \\
&=&\mathfrak{x}^{n}+\tsum\nolimits_{j=1}^{(n-1)/2}\left( -1\right) ^{j}%
\tfrac{n}{n-j}\tbinom{n-j}{j}\mathfrak{x}^{n-2j}.\text{ \ }\blacksquare
\end{eqnarray*}

\bigskip

We have proved:

\begin{theorem}
\label{presentation}Let $g$ be a positive integer, $p$ a prime, $q=p^{g}$
and set $\mathfrak{G}=SL_{2}(\mathbb{F}_{q})$. Denote by $\mathfrak{X}$ the
standard representation of $\mathfrak{G}$ on $\mathbb{F}_{q}^{2}$ and let $%
\mathfrak{x}$ be an indeterminate over $%
\mathbb{Z}
$. The assignment $\mathfrak{X\longmapsto x}$\ induces an isomorphism of
rings:%
\begin{equation*}
K_{0}\left( \mathfrak{G}\right) \simeq \frac{%
\mathbb{Z}
\lbrack \mathfrak{x}]}{\left( \mathfrak{f}^{[g]}(\mathfrak{x})-\mathfrak{x}%
\right) 
\mathbb{Z}
\lbrack \mathfrak{x}]},
\end{equation*}

\noindent where $\mathfrak{f}^{[g]}(\mathfrak{x})=\left( \mathfrak{f}\circ 
\mathfrak{f}\circ ...\circ \mathfrak{f}\right) \left( \mathfrak{x}\right) $
is the monic polynomial of $%
\mathbb{Z}
\lbrack \mathfrak{x}]$ having degree $p^{g}$ that is obtained by composing $%
g $-times wit itself the monic degree $p$ polynomial:%
\begin{equation*}
\mathfrak{f}(\mathfrak{x}):=\dsum\nolimits_{j=0}^{\left\lfloor
p/2\right\rfloor }\left( -1\right) ^{j}\dfrac{p}{p-j}\dbinom{p-j}{j}%
\mathfrak{x}^{p-2j}.
\end{equation*}
\end{theorem}

\bigskip

At the time of writing of this paper, we do not know much about the
properties of the polynomial $\mathfrak{f}^{[g]}(\mathfrak{x})-\mathfrak{x}$
when viewed over $%
\mathbb{Z}
$. Notice that if $p>2$, $\mathfrak{f}^{[1]}(\mathfrak{x})-\mathfrak{x}$ is
an odd polynomial; using the easy to check facts that for any integer $n\geq
0$ we have $\mathfrak{m}_{n}\left( 2\right) =n+1,$ and that:%
\begin{equation*}
\mathfrak{m}_{n}\left( 1\right) =\left\{ 
\begin{array}{c}
-1\text{, if }n\equiv 3,4(\func{mod}6) \\ 
0\text{, \ if }n\equiv 2,5(\func{mod}6) \\ 
1\text{, if }n\equiv 0,1(\func{mod}6),%
\end{array}%
\right.
\end{equation*}

\noindent one deduce that $\mathfrak{f}^{[1]}(\mathfrak{x})-\mathfrak{x}$
always admits $0,\pm 1,\pm 2$ as roots, as long as $p>3$ (roots $0$ and $\pm
2$ also occur for $p=3$). Furthermore, from computer elaborations, $%
\mathfrak{f}^{[1]}(\mathfrak{x})-\mathfrak{x}$ seems to have only real roots.

In general, it is natural to ask what we can say about the irreducible
factors over $%
\mathbb{Q}
$ of $\mathfrak{f}^{[g]}(\mathfrak{x})-\mathfrak{x}$. We do not have an
answer for this. Nevertheless, after tensoring $K_{0}\left( \mathfrak{G}%
\right) $ with $%
\mathbb{Z}
_{p}$, we can prove:

\begin{corollary}
\noindent \label{cor_presentation}Let $K_{0}\left( \mathfrak{G}\right)
_{p}=K_{0}\left( \mathfrak{G}\right) \otimes _{%
\mathbb{Z}
}%
\mathbb{Z}
_{p}$. For any positive divisor $d$ of $g$, let $\psi (d)$ be the number of
monic irreducible polynomials of degree $d$ in $\mathbb{F}_{p}[\mathfrak{x}]$%
. Then:

\begin{description}
\item[(a)] The special fiber of $K_{0}\left( \mathfrak{G}\right) _{p}$ is a
split$\mathbb{\ F}_{p}$-algebra isomorphic to $\dprod\nolimits_{d|g}\mathbb{F%
}_{p^{d}}^{\psi (d)};$

\item[(b)] The generic fiber of $K_{0}\left( \mathfrak{G}\right) _{p}$ is a
split $%
\mathbb{Q}
_{p}$-algebra isomorphic to $\dprod\nolimits_{d|g}\mathbb{%
\mathbb{Q}
}_{p^{d}}^{\psi (d)}.$
\end{description}

\noindent (Here we denoted by $\mathbb{%
\mathbb{Q}
}_{p^{d}}$ the degree $d$ unramified extension of $\mathbb{%
\mathbb{Q}
}_{p}$ contained in a fixed algebraic closure of $\mathbb{%
\mathbb{Q}
}_{p}$).
\end{corollary}

\textbf{Proof }By the explicit formula given above for $\mathfrak{f}(%
\mathfrak{x})$, we see that $\mathfrak{f}(\mathfrak{x})\equiv \mathfrak{x}%
^{p}(\func{mod}p%
\mathbb{Z}
\lbrack \mathfrak{x}])$: this is clear if $p=2$, otherwise notice that $%
\tfrac{p}{p-j}\tbinom{p-j}{j}=p\cdot \tfrac{(p-j-1)!}{j!\cdot (p-2j)!}$ and
the last denominator is prime to $p$ if $1\leq j\leq \tfrac{p-1}{2}$,
implying that $\tfrac{(p-j-1)!}{j!\cdot (p-2j)!}\in 
\mathbb{Z}
$. We conclude that 
\begin{equation*}
\mathfrak{f}^{[g]}(\mathfrak{x})-\mathfrak{x\equiv x}^{q}-\mathfrak{x}(\func{%
mod}p%
\mathbb{Z}
\lbrack \mathfrak{x}])
\end{equation*}
and statement (a) follows. Part (b)\ follows from (a) and Hensel's lemma. $%
\blacksquare $

\bigskip

\begin{remark}
We also have isomorphisms of algebras: $K_{0}\left( \mathfrak{G}\right)
_{p}\otimes _{%
\mathbb{Z}
_{p}}\mathbb{F}_{q}\simeq \left( \mathbb{F}_{q}\right) ^{q}$ and $%
K_{0}\left( \mathfrak{G}\right) _{p}\otimes _{%
\mathbb{Z}
_{p}}\mathbb{%
\mathbb{Q}
}_{q}\simeq \left( \mathbb{%
\mathbb{Q}
}_{q}\right) ^{q}$.
\end{remark}

It is interesting to notice that we can give an explicit formula also for $%
\mathfrak{f}^{[g]}(\mathfrak{x})$. As the following proposition uses Serre's
relation $(\Sigma _{g})$, it seems that an explicit formula for $\mathfrak{f}%
^{[i]}(\mathfrak{x})$ when $i\neq 1,g$ would probably require more work.

\begin{proposition}
\label{gth compo}We have $\mathfrak{f}^{[g]}(\mathfrak{x})=\mathfrak{m}_{q}(%
\mathfrak{x})-\mathfrak{m}_{q-2}(\mathfrak{x})$, so that:%
\begin{equation*}
\mathfrak{f}^{[g]}(\mathfrak{x})=\dsum\nolimits_{j=0}^{\left\lfloor
q/2\right\rfloor }\left( -1\right) ^{j}\dfrac{q}{q-j}\dbinom{q-j}{j}%
\mathfrak{x}^{q-2j}.
\end{equation*}
\end{proposition}

\textbf{Proof }Let $\tilde{\pi}:%
\mathbb{Z}
\lbrack \mathfrak{x}]\rightarrow K_{0}\left( \mathfrak{G}\right) $ be the
epimorphism of rings obtained by sending $\mathfrak{x}$ to $\mathfrak{X}$.
Relation $(\Sigma _{g})$ implies that $\mathfrak{M}_{1}=\mathfrak{M}_{q}-%
\mathfrak{M}_{q-2}$ in $K_{0}(\mathfrak{G})$, that is $\mathfrak{M}_{q}-%
\mathfrak{M}_{q-2}-\mathfrak{X}=0$. This means, by formula (\ref{meaningful1}%
), that $\mathfrak{X}$ satisfies the polynomial $\mathfrak{m}_{q}(\mathfrak{x%
})-\mathfrak{m}_{q-2}(\mathfrak{x})-\mathfrak{x\in 
\mathbb{Z}
}[\mathfrak{x}]$, so that $\mathfrak{m}_{q}(\mathfrak{x})-\mathfrak{m}_{q-2}(%
\mathfrak{x})-\mathfrak{x\in }\ker \tilde{\pi}=\left( \mathfrak{f}^{[g]}(%
\mathfrak{x})-\mathfrak{x}\right) 
\mathbb{Z}
\lbrack \mathfrak{x}]$. Since $\mathfrak{m}_{q}(\mathfrak{x})-\mathfrak{m}%
_{q-2}(\mathfrak{x})-\mathfrak{x}$ and $\mathfrak{f}^{[g]}(\mathfrak{x})-%
\mathfrak{x}$ are both monic of degree $q$, the last relation implies that
they have to be equal and $\mathfrak{f}^{[g]}(\mathfrak{x})=\mathfrak{m}_{q}(%
\mathfrak{x})-\mathfrak{m}_{q-2}(\mathfrak{x})$. The proposition now follows
from Corollary \ref{n-n2}. $\blacksquare $

\section{A conjecture}

The following fact was pointed out to us by G. Savin:

\begin{proposition}
\label{tensorp}Let $p$ be a prime and $q=p^{g}>1$ be an integral power of $p$%
. Let $\mathbb{G}$ be a simply connected semisimple algebraic group defined
and split over $\mathbb{F}_{q}$, and denote by $K_{0}(\mathbb{G)}$ the
Grothendieck ring of $\mathbb{F}_{q}[\mathbb{G]}$-rational modules of finite 
$\mathbb{F}_{q}$-dimension. If $\mathfrak{M}$ is an element of $K_{0}(%
\mathbb{G)}$ and $i$ is any non-negative integer, we have:%
\begin{equation*}
\mathfrak{M}^{[i]}\equiv \mathfrak{M}^{p^{i}}(\func{mod}pK_{0}(\mathbb{G))}.
\end{equation*}
\end{proposition}

\textbf{Proof }Let $\mathbb{T}$\ be a maximal torus of $\mathbb{G}$ defined
and split over $\mathbb{F}_{q}$,\ and denote by $X=X(\mathbb{T})$ its
character group. For any $\lambda \in X$, denote by $e(\lambda )$ the
corresponding basis element of the group ring $%
\mathbb{Z}
\lbrack X]$, so that $e(\lambda +\lambda ^{\prime })=e(\lambda )e(\lambda
^{\prime })$ for any characters $\lambda $ and $\lambda ^{\prime }$.\ 

Fix a $\mathbb{G}$-module $\mathfrak{M}$ and write its formal character as:%
\begin{equation*}
\limfunc{ch}\mathfrak{M}=\dsum\nolimits_{\lambda \in X}m_{\lambda }\cdot
e(\lambda ),
\end{equation*}

\noindent where $m_{\lambda }$ is the dimension of the $\lambda $-isotypic
submodule of $\mathfrak{M}$. For a positive integer $i$, the $p^{i}$th power
automorphism of $\mathbb{\bar{F}}_{q}$ induces an action on $%
\mathbb{Z}
\lbrack X]$ by sending a basis element $e(\lambda )$ to $e(p^{i}\lambda )$,
so that: 
\begin{eqnarray*}
\limfunc{ch}(\mathfrak{M}^{[i]}) &=&\dsum\nolimits_{\lambda \in X}m_{\lambda
}\cdot e(\lambda )^{p^{i}} \\
&\equiv &\left( \dsum\nolimits_{\lambda \in X}m_{\lambda }\cdot e(\lambda
)\right) ^{p^{i}}(\func{mod}p%
\mathbb{Z}
\lbrack X]\mathbb{)}\text{.}
\end{eqnarray*}

The formal character $\left( \tsum\nolimits_{\lambda \in X}m_{\lambda }\cdot
e(\lambda )\right) ^{p^{i}}$ is the element associated to $\mathfrak{M}%
^{p^{i}}$ by the map $\limfunc{ch}:K_{0}(\mathbb{G})\longrightarrow 
\mathbb{Z}
\lbrack X]$. We have therefore:%
\begin{equation}
\limfunc{ch}(\mathfrak{M}^{[i]})\equiv \limfunc{ch}(\mathfrak{M}^{p^{i}})%
\text{ \ }(\func{mod}p%
\mathbb{Z}
\lbrack X]\mathbb{)}\text{.}  \tag{$1$}  \label{11}
\end{equation}

Let $\mathcal{W}$ denotes the Weyl group of the pair $(\mathbb{G},\mathbb{T}%
) $. By \cite{Jan}\ II.5.8, the map $\limfunc{ch}$ induces an isomorphism of
commutative rings:%
\begin{equation*}
\limfunc{ch}:K_{0}(\mathbb{G})\overset{\sim }{\longrightarrow }%
\mathbb{Z}
\lbrack X]^{\mathcal{W}}.
\end{equation*}

Write $\limfunc{ch}(\mathfrak{M}^{[i]})=\tsum\nolimits_{\lambda \in
X}a_{\lambda }\cdot e(\lambda )$ and $\limfunc{ch}(\mathfrak{M}%
^{p^{i}})=\tsum\nolimits_{\lambda \in X}b_{\lambda }\cdot e(\lambda )$, so
that $\limfunc{ch}(\mathfrak{M}^{[i]}-\mathfrak{M}^{p^{i}})=\tsum\nolimits_{%
\lambda \in X}\left( a_{\lambda }-b_{\lambda }\right) \cdot e(\lambda )$ is
such that:%
\begin{equation}
\dsum\nolimits_{\lambda \in X}\left( a_{\lambda }-b_{\lambda }\right) \cdot
e(\lambda )=\dsum\nolimits_{\lambda \in X}\left( a_{\lambda }-b_{\lambda
}\right) \cdot e(w\lambda )  \tag{$2$}  \label{22}
\end{equation}%
for all $w\in \mathcal{W}$.

By (\ref{11}), there are integers $c_{\lambda }$ such that $a_{\lambda
}-b_{\lambda }=pc_{\lambda }$ for all $\lambda \in X$. Since $%
\mathbb{Z}
\lbrack X]$ is $%
\mathbb{Z}
$-flat, we can view it as a subring of $%
\mathbb{Q}
[X]$, in which we have, for any $w\in \mathcal{W}$:%
\begin{eqnarray*}
w\cdot \left( \dsum\nolimits_{\lambda \in X}c_{\lambda }\cdot e(\lambda
)\right) &=&\frac{1}{p}\dsum\nolimits_{\lambda \in X}\left( a_{\lambda
}-b_{\lambda }\right) \cdot e(w\lambda ) \\
&=&\frac{1}{p}\dsum\nolimits_{\lambda \in X}\left( a_{\lambda }-b_{\lambda
}\right) \cdot e(\lambda ),
\end{eqnarray*}

\noindent where the last equality follows from (\ref{22}). Therefore $w\cdot
\left( \tsum\nolimits_{\lambda \in X}c_{\lambda }\cdot e(\lambda )\right)
=\tsum\nolimits_{\lambda \in X}c_{\lambda }\cdot e(\lambda )$ in $%
\mathbb{Z}
\lbrack X]$\ for all $w\in \mathcal{W}$ and 
\begin{equation*}
\limfunc{ch}(\mathfrak{M}^{[i]}-\mathfrak{M}^{p^{i}})\in p%
\mathbb{Z}
\lbrack X]^{\mathcal{W}}.
\end{equation*}

This\ implies that $\mathfrak{M}^{[i]}\ $is congruent to $\mathfrak{M}%
^{p^{i}}$ modulo the ideal generated by $p$ in $K_{0}(\mathbb{G)}$. $%
\blacksquare $

\bigskip

Motivated by Theorem \ref{presentation}, Corollary \ref{cor_presentation}
and Proposition \ref{tensorp}, we are led to the following:

\begin{conjecture}
\label{conJ}Let $p$ be a prime and $q=p^{g}>1$ be an integral power of $p$.
Let $\mathbb{G}$ be a simply connected semisimple algebraic group defined
and split over $\mathbb{F}_{q}$, whose rank is $\ell >0$. Denote by $K_{0}(%
\mathbb{G}(\mathbb{F}_{q}))$ the Grothendieck ring of finitely generated $%
\mathbb{F}_{q}[\mathbb{G}(\mathbb{F}_{q})]$-modules. Then there exist $\ell $
monic polynomials $\mathfrak{f}_{1}(\mathfrak{x}_{1})\in 
\mathbb{Z}
\lbrack \mathfrak{x}_{1}],...,\mathfrak{f}_{\ell }(\mathfrak{x}_{\ell })\in 
\mathbb{Z}
\lbrack \mathfrak{x}_{\ell }]$ having degree $p$\ such that:%
\begin{equation*}
K_{0}\left( \mathbb{G}(\mathbb{F}_{q})\right) \simeq \frac{%
\mathbb{Z}
\lbrack \mathfrak{x}_{1},...,\mathfrak{x}_{\ell }]}{\left( \mathfrak{f}%
_{1}^{[g]}(\mathfrak{x}_{1})-\mathfrak{x}_{1}\mathfrak{,...,f}_{\ell }^{[g]}(%
\mathfrak{x}_{\ell })-\mathfrak{x}_{\ell }\right) 
\mathbb{Z}
\lbrack \mathfrak{x}_{1},...,\mathfrak{x}_{\ell }]},
\end{equation*}

\noindent where for any $i$, $1\leq i\leq \ell $, $\mathfrak{f}_{i}^{[g]}(%
\mathfrak{x}_{i})$ is the polynomial obtained by composing $\mathfrak{f}_{i}(%
\mathfrak{x}_{i})$ with itself $g$ times.

Furthermore, $\mathfrak{f}_{i}(\mathfrak{x}_{i})\equiv \mathfrak{x}_{i}^{p}(%
\func{mod}p%
\mathbb{Z}
\lbrack \mathfrak{x}_{i}])$ for any $i$, $1\leq i\leq \ell $.
\end{conjecture}

\bigskip

The idea behind the above statement is that if $\pi $ is an isomorphism of $%
\mathbb{Z}
\lbrack \mathfrak{x}_{1},...,\mathfrak{x}_{\ell }]/(\mathfrak{f}_{1}^{[g]}(%
\mathfrak{x}_{1})-\mathfrak{x}_{1}\mathfrak{,...,f}_{\ell }^{[g]}(\mathfrak{x%
}_{\ell })-\mathfrak{x}_{\ell })$ onto $K_{0}\left( \mathbb{G}(\mathbb{F}%
_{q})\right) $, and if we set $\mathfrak{X}_{i}:=\pi (\mathfrak{x}_{i})$ for 
$1\leq i\leq \ell $, then $\mathfrak{f}_{i}(\mathfrak{X}_{i})\in K_{0}\left( 
\mathbb{G}(\mathbb{F}_{q})\right) $ should be the Frobenius twist $\mathfrak{%
X}_{i}^{[1]}$. This means that the relations imposed above in the algebra $%
\mathbb{Z}
\lbrack \mathfrak{x}_{1},...,\mathfrak{x}_{\ell }]$ are the obvious ones
that translate into $\mathfrak{X}_{i}^{[g]}=\mathfrak{X}_{i}$ for all $i$.

Here is some evidence for the conjecture:

\begin{description}
\item[(a)] As proved in the previous paragraph, the conjecture is true for $%
\mathbb{G}=SL_{2}$ over $\mathbb{F}_{q}$ ($\ell =1$), in which case we can
also give an explicit formula for the polynomial $\mathfrak{f}(\mathfrak{x})=%
\mathfrak{f}_{1}(\mathfrak{x}_{1})$ (Theorem \ref{presentation}).

\item[(b)] A theorem of Steinberg (\cite{St}) states that if $\mathbb{G}$ is
a simply connected semisimple algebraic group over $\mathbb{F}_{q}$, the
number of semisimple conjugacy classes of $\mathbb{G}(\mathbb{F}_{q})$ is
equal to $q^{\ell }$, where $\ell $ is the rank of $\mathbb{G}$. Therefore $%
K_{0}\left( \mathbb{G}(\mathbb{F}_{q})\right) \simeq 
\mathbb{Z}
^{q^{\ell }}$ as $%
\mathbb{Z}
$-modules, which follows from the conjecture.

\item[(c)] Since $\mathfrak{h}(\mathfrak{x}_{1},...,\mathfrak{x}_{\ell
})^{q}=\mathfrak{h(}\mathfrak{x}_{1}^{q},...,\mathfrak{x}_{\ell }^{q}%
\mathfrak{)}$ for any polynomial $\mathfrak{h}(\mathfrak{x}_{1},...,%
\mathfrak{x}_{\ell })\in \mathbb{F}_{q}[\mathfrak{x}_{1},...,\mathfrak{x}%
_{\ell }]$, the conjecture implies that $\overline{\mathfrak{M}}^{q}=%
\overline{\mathfrak{M}}$ for any $\overline{\mathfrak{M}}\in K_{0}(\mathbb{G}%
(\mathbb{F}_{q}))\otimes _{%
\mathbb{Z}
}\mathbb{F}_{q}$. This fact is also a consequence of Proposition \ref%
{tensorp}.

\item[(d)] Assume we are given a surjective homomorphism of $\mathbb{F}_{q}$%
-algebras:%
\begin{equation*}
\gamma :\mathbb{F}_{q}[\mathfrak{x}_{1},...,\mathfrak{x}_{\ell
}]\longrightarrow K_{0}(\mathbb{G}(\mathbb{F}_{q}))\otimes _{%
\mathbb{Z}
}\mathbb{F}_{q}\text{.}
\end{equation*}%
Proposition \ref{tensorp} implies that $\gamma (\mathfrak{x}_{i})^{q}=\gamma
(\mathfrak{x}_{i})$ for any integer $i$, $1\leq i\leq \ell $; in particular
the kernel of $\gamma $ contains the ideal generated by the polynomials $%
\mathfrak{x}_{1}^{q}-\mathfrak{x}_{1}\mathfrak{,...,x}_{\ell }^{q}-\mathfrak{%
x}_{\ell }$. By dimension reasons we must have an isomorphism of $\mathbb{F}%
_{q}$-algebras:%
\begin{equation*}
\frac{\mathbb{F}_{q}[\mathfrak{x}_{1},...,\mathfrak{x}_{\ell }]}{\left( 
\mathfrak{x}_{1}^{q}-\mathfrak{x}_{1}\mathfrak{,...,x}_{\ell }^{q}-\mathfrak{%
x}_{\ell }\right) \mathbb{F}_{q}[\mathfrak{x}_{1},...,\mathfrak{x}_{\ell }]}%
\overset{\sim }{\longrightarrow }K_{0}(\mathbb{G}(\mathbb{F}_{q}))\otimes _{%
\mathbb{Z}
}\mathbb{F}_{q}.
\end{equation*}%
This is predicted by Conjecture \ref{conJ}.
\end{description}

If Conjecture \ref{conJ} is correct, one would like to determine explicit
formulae for the polynomials $\mathfrak{f}_{1}(\mathfrak{x}_{1}),...,%
\mathfrak{f}_{\ell }(\mathfrak{x}_{\ell })$ and to relate factorization
properties of these polynomials in $%
\mathbb{Z}
\lbrack \mathfrak{x}_{i}]$ to algebraic properties of the group $\mathbb{G}(%
\mathbb{F}_{q})$.

\bibliographystyle{amsplain}
\bibliography{davide2}
\bigskip

\end{document}